\documentclass[twoside,a4paper,leqno,12pt]{article}
\usepackage{amsmath,amsthm,amsfonts,amssymb}

\makeatletter
\def\@seccntformat#1{\csname the#1\endcsname.\quad}
\renewcommand\section{\@startsection {section}{1}{\z@}%
                                   {-3.5ex \@plus -1ex \@minus -.2ex}%
                                   {2.3ex \@plus.2ex}%
                                    {\normalfont\bf\flushleft }}
\makeatother

\date{\today}

\newcommand{\R}{\mathbf{R}}

\newcommand{\C}{\mathbf{C}}  
  
\newcommand{\Hy}{\mathbf{H}}

\newcommand{\F}{\mathcal{F}}
\newcommand{\Bo}{\mathcal{B}}
\newcommand{\dd}{\mathrm{d}}

\newcommand{\pa}{\partial}

\newcommand{\eil}{\overset{(\text{\rm law})}{=}}
\newcommand{\wt}{\widetilde}
\newcommand{\wh}{\widehat}
\theoremstyle{plain}
  \newtheorem{thm}{Theorem}[section]
  \newtheorem{cor}[thm]{Corollary}
  
  \newtheorem{prop}[thm]{Proposition}

\theoremstyle{remark}
  \newtheorem{rem}{Remark}[section]
\numberwithin{equation}{section}

\allowdisplaybreaks[4]

\begin{document}

\vspace*{2cm} 

\begin{center}
{\Large\bf Limiting behaviors of the Brownian motions 
on hyperbolic spaces}
\end{center}

\bigskip

\begin{center}
H. MATSUMOTO (Nagoya) 
\footnote{This research was partially supported by 
Grant-in-Aid for Scientific Research (No. 19204010), 
Japan Society for the Promotion of Science.}
\end{center}

\bigskip

\begin{abstract}
Using the explicit representations of the Brownian motions 
on the hyperbolic spaces, 
we show that their almost sure convergence 
and the central limit theorems for the radial components 
as time tends to infinity are easily obtained.    
We also give a straightforward strategy to obtain 
the explicit expressions for the limit distributions 
or the Poisson kernels.
\end{abstract} 

\bigskip

\section{Introduction}

Hyperbolic spaces are non-compact Riemannian symmetric spaces 
of rank one.   
By classification, we have four types of hyperbolic spaces: 
the real one $\Hy_r^n=SO_0(1,n)/SO(n)$, 
the complex one $\Hy_c^n=SU(1,n)/SU(n)$, 
the quaternionic one $Sp(1,n)/(Sp(1)\times Sp(n))$ 
and the Cayley hyperbolic plane.   
In this article we consider the limiting behaviors 
of the Brownian motions, that is, the diffusion processes 
generated by the Laplace-Beltrami operators 
on the first three types of the hyperbolic spaces.   

The hyperbolic spaces have negative bounded curvatures.   
The Brownian motions on negatively curved manifolds have been studied 
in the connection of so-called the Liouville property 
by many authors and it is well known that 
the Brownian motions tends to infinity almost surely 
as time tends to infinity. 
See, e.g., Kifer \cite{kifer}.  
Needless to say, 
the limit distributions are given by the Poisson kernels.   

On the other hand, 
the Brownian motions on the Riemannian symmetric spaces 
of non-compact type have been also studied by several authors 
since the work by Malliavin-Malliavin \cite{mm}.   
Among them, we refer to Babillot \cite{babillot}, 
where a central limit theorem for the radial components 
of the Brownian motions has been shown.  

The purpose of this article is to show that 
these properties are easily and directly shown 
for the Brownian motions on the hyperbolic spaces 
if we adopt the upper half space realizations 
of the hyperbolic spaces instead of the ball models.   
We can describe the same stories on the three types of spaces.   
At first, 
by solving the corresponding stochastic differential equations, 
we represent the Brownian motions in closed forms 
as Wiener functionals.   
Then, the almost sure convergence of them is readily seen 
from the representations.   
Moreover, by inserting the representations 
into the formulae for the distance functions, 
we can also show 
the central limit theorems for the radial components.   

For the computations of the limiting distributions 
or the Poisson kernels, 
we need some results on the distributions of the random variables 
defined by the perpetual (infinite) integrals in time 
of the usual geometric Brownian motions with negative drifts.   
The auxiliary results are given in the appendix and, by using them, 
we compute the Fourier transforms of the limiting distributions 
and the inverse transforms in direct ways.   
\section{Real hyperbolic spaces}

For $n\geqq 1$, let $\Hy_r^{n+1}$ be the upper half space 
in $\R^{n+1}$, 
\begin{equation*}
\Hy_r^{n+1} = \{z=(x,y)=(x_1,...,x_n,y); x\in\R^n, y>0\}, 
\end{equation*}
endowed with the Riemannian metric $\dd s^2=y^{-2}(\dd x^2+\dd y^2)$. 
The volume element is given by $y^{-n-1}\dd x \dd y$, 
the distance function $d(z,z')$ is given by 
\begin{equation} \label{re:dist}
\cosh(d(z,z')) = \frac{|x-x'|^2+y^2+(y')^2}{2yy'},
\end{equation}
in an obvious notation, where $|x|$ is the Euclidean norm.  
The Laplace-Beltrami operator is written as 
\begin{equation*}
\Delta_r = y^2 \sum_{j=1}^n \frac{\pa^2}{\pa x_j^2} + 
y^2 \frac{\pa^2}{\pa y^2} - (n-1)y \frac{\pa}{\pa y}.
\end{equation*}
For details of the fundamental objects on $\Hy_r^{n+1}$, 
see, e.g., Davies \cite{davies}.   

We first show an explicit expression as Wiener functional 
of the Brownian motion on $\Hy_r^{n+1}$ 
by soloving the corresponding stochastic differential equation.   
Let $(W^{(n+1)}, \Bo^{(n+1)}, P^{(n+1)})$ be 
the $(n+1)$-dimensional standard Wiener space 
with the canonical filtration $\{\Bo_s^{(n+1)}\}_{t\geqq0}$.   
Corresponding to the rectangular coordinate, 
we denote an element of $W^{(n+1)}$ by 
\begin{equation*}
(w(\cdot),B(\cdot)) \qquad \text{\rm or} \qquad 
(w_1(\cdot),...,w_n(\cdot),B(\cdot)), 
\end{equation*}
which is an $\R^{n+1}$-valued continuous function 
on $[0,\infty)$ with $w_i(0)=B(0)=0$.   
Then the Brownian motion on $\Hy_r^{n+1}$, 
the diffusion process with infinitesimal generator $\Delta_r/2$, 
is obtained by solving the following stochastic differential equation 
defined on $(W^{(n+1)}, \Bo^{(n+1)}, P^{(n+1)})$ (see \cite{iw}): 
\begin{align*}
 & \dd X_i(t) = Y(t) \dd w_i(t), \qquad \qquad \qquad \qquad 
i=1,...,n, \\
 & \dd Y(t) = Y(t) \dd B(t) - \frac{n-1}{2} Y(t) \dd t.
\end{align*}
The unique solution 
$Z_z=\bigl\{\bigl(X(t,z),Y(t,z)\bigr)\bigr\}_{t\geqq0}, z=(x,y),$ 
satisfying $X(0)=x$ and $Y(0)=y$ is given by 
\begin{align*}
 & X_i(t,z) = x_i + \int_0^t y \exp(B_s^{(-\mu)}) \dd w_i(s), \\
 & Y(t,z) = y \exp(B_t^{(-\mu)}),
\end{align*}
where $B_s^{(-\mu)}=B(s)-\mu s$ and 
$\mu=n/2$.  
$\{Y(t,z)\}$ is a usual geometric Brownian motion with negative drift 
and it is easy to see that $Z_z(t)$ converges to the boundary 
as $t\to\infty$ almost surely.   

Now we consider the exponential functional $A_t^{(-\mu)}$ given by 
\begin{equation*}
A_t^{(-\mu)} = \int_0^t\exp(2B_s^{(-\mu)})\dd s.  
\end{equation*}
Then it is easy to see the identity in law 
\begin{equation*}
(X(t,z),Y(t,z)) \eil \bigl(x+yw\bigl(A_t^{(-\mu)}\bigr), 
y\exp\bigl(B_t^{(-\mu)}\bigr) \bigr)
\end{equation*}
for fixed $t>0$.   

An explicit expression for the density of 
the distribution of $(A_t^{(-\mu)},B_t^{(-\mu)})$ is 
known by Yor \cite{yor-some} and, by using it, 
Gruet \cite{gruet} has shown an expression 
for the heat kernel of the semigroup 
generated by $\Delta_r$.   
For the classical expression. see Davies \cite{davies}.   
We also refer to \cite{ag, im-jfa99, m, survey2} 
for related topics.   

We combine the identity in law with formula \eqref{re:dist}.   
Then we get 
\begin{equation*}
 \cosh(d(Z(t,z),z)) \eil 
\frac{1}{2} \biggl\{|yw\bigl(A_t^{(-\mu)}\bigr)|^2+1\biggr\} 
\exp\bigl(-B_t^{(-\mu)}\bigr) + 
\frac{1}{2} \exp\bigl(B_t^{(-\mu)}\bigr).
\end{equation*}
Since $A_t^{(-\mu)}$ converges as $t\to\infty$ and 
$\log(\cosh(u))=u\cdot(1+o(1))$ as $u\to\infty$, 
we readily get the following central limit theorem.   

\medskip

\begin{thm} \label{rt:clt}
The probability distribution of $\sqrt{t}^{-1}(d(Z(t,z),z)-nt/2)$ 
converges weakly as $t\to\infty$ to the standard normal distribution.
\end{thm}

\medskip

Recall formula $\Delta_r d(z_0,\cdot)=n\coth d(z_0,\cdot)$.  
Then, by It\^o's formula, we get 
\begin{equation*} \begin{split}
d(Z_z(t),z) = & \sum_{i=1}^n \int_0^t \frac{1}{\sinh d(z,Z_z(s))} 
\frac{X_z^i(s)-x}{y} \dd w_s^i  \\
 & + \int_0^t \frac{1}{\sinh d(z,Z_z(s))} \biggl( \frac{Y_z(s)}{y} - 
\cosh d(z,Z_z(s)) \biggr) \dd B(s) \\
 & + \frac{n}{2} \int \coth d(z,Z_z(s)) \dd s,
\end{split} \end{equation*}
from which the theorem may also be proven.   

Next we recall 
Dufresne's identity (Theorem \ref{pt:marg} in the Appendix) 
in law $A_\infty^{(-\mu)} \eil (2\gamma_\mu)^{-1}$ for a 
Gamma random variable $\gamma_\mu$ with parameter $\mu.$  
Then, for a bounded continuous function $\varphi$ on $\R^n$, 
we obtain 
\begin{align*}
 & E[\varphi(X(t,z))] = 
E[\varphi\bigl(x+yw\bigl(A_t^{(-n/2)}\bigr)\bigr)] \\
 & \to \int_0^\infty \frac{1}{\Gamma(n/2)} t^{(n/2)-1} e^{-t}\dd t 
\int_{\R^n} \varphi(x+\eta) \frac{1}{(2\pi y^2/2t)^{n/2}} 
\exp\biggl(-\frac{|\eta|^2}{2y^2/2t}\biggr) \dd\eta \\
 & = \int_{\R^n} \varphi(\xi)\dd\xi 
 \int_0^\infty \frac{1}{\Gamma(n/2)}
 \frac{1}{\pi^{n/2}y^n} t^{n-1} 
 \exp\biggl(-\frac{y^2+|\xi-x|^2}{y^2}t\biggr) \dd t \\
 & = \int_{\R^n} \varphi(\xi) p_{n+1}(\xi-x,y) \dd\xi,
\end{align*}
where 
\begin{equation*}
p_{n+1}(\xi,y) = \frac{2^{n-1}\Gamma((n+1)/2)}{\pi^{(n+1)/2}} 
\frac{y^n}{(y^2+|\xi|^2)^n}, \qquad \xi\in\R^n, 
\end{equation*}
and we have used the duplication formula for the Gamma function.  

Hence we have proved the following.  

\medskip

\begin{thm} \label{rt:poisson}
For any $(x,y)\in\Hy_r^{n+1},$ $X(t,z)$ converges almost surely 
as $t\to\infty$ and the density of the limit distribution is 
given by the Poisson kernel $p_{n+1}(\xi-x,y).$  
In particular, when $n=1,$ 
the limit distribution is Cauchy{\rm .}
\end{thm}

\medskip

We end this section with mentioning on the Poisson kernel 
in the Euclidean spaces and on the Fourier transforms.   
The Poisson kernel on $\R^{n+1}$ of the hyperplane $\{y=0\}$ 
is given by 
\begin{equation*}
q_{n+1}(\xi,y) = \frac{\Gamma((n+1)/2)}{\pi^{(n+1)/2}} 
\frac{y}{(y^2+|\xi|^2)^{(n+1)/2}} , 
\end{equation*}
which is different from $p_{n+1}(\xi,y)$ for $n\geqq2$.   
The Brownian motion on the hyperbolic plane $\Hy^2$ is 
a time change of the $2$-dimensional standard Brownian motion, 
and the Poisson kernels coincide.  

It is well known that the Fourier transform of $q_{n+1}(\xi,y)$ 
in $\xi$ is the simple exponential function, 
\begin{equation*}
\int_{\R^n} e^{\sqrt{-1}\langle \lambda, \xi \rangle} 
q_{n+1}(\xi,y) \dd\xi = e^{-y|\lambda|}. 
\end{equation*}
For the hyperbolic spaces, we can show, for example, 
\begin{align*}
 & \varphi_3(\lambda;y) \equiv 
\int_{\R^2} e^{\sqrt{-1}\langle \lambda, \xi \rangle} 
p_3(\xi,y)\dd\xi = y|\lambda|K_1(y|\lambda|) \\
\intertext{and}
 & \varphi_4(\lambda;y) \equiv 
\int_{\R^3} e^{\sqrt{-1}\langle \lambda, \xi \rangle} 
p_4(\xi,y)\dd\xi = (y|\lambda|+1) e^{-y|\lambda|},
\end{align*}
where $K_1$ is the modified Bessel function.   
By virtue of the strong Markov property, 
we can easily show that the distribution of $X(\tau_a)$ 
for the first hitting time $\tau_a$ 
of the Brownian motion $\{Z_z(t,z)\}$ at the level $y=a, a>0$ 
is determined by the characteristic function written by 
\begin{equation*}
E[\exp(\sqrt{-1} \langle \lambda, X(\tau_a) \rangle )] = 
e^{\sqrt{-1}\langle \lambda,x \rangle} 
\frac{\varphi_n(\lambda;y)}{\varphi_n(\lambda;a)}, \qquad 
\lambda \in \R^n. 
\end{equation*}
\section{Complex hyperbolic spaces}

Let $\Hy_c^n, n\geqq 2,$ be the upper half space of $\C^n$ given by 
\begin{equation*}
\{z=(z_1,z_2,...,z_n)=(z_1,\widetilde{z})\in\C^n; 
h(z)\equiv \mathrm{Im}(z_1)-|\widetilde{z}|^2>0\}, 
\end{equation*}
endowed with the Bergmann metric
\begin{equation*}
\dd s^2 = - \sum_{j,k=1}^n \pa_{z_j} \pa_{\overline{z_k}} 
\log(h) \dd z_j \dd \overline{z_k}.
\end{equation*}
The unit ball $\{|z|<1\}$ in $\C^n$ with the Bergmann metric 
\begin{equation*}
- \sum_{j,k=1}^n \pa_{z_j} \pa_{\overline{z_k}} 
\log(1-|z|^2) \dd z_j \dd \overline{z_k}.
\end{equation*}
is isometric with $\Hy_c^n$.   
For details, we refer to \cite{deb1,deb2,goldman,venkov}. 
We should be aware of difference of conventions.   
The curvatures of these manifolds are bounded and negative, 
but are not constant (cf. p.190, \cite{goldman}).  

We change the first coordinate 
by $x_1=\mathrm{Re}(z_1)/2$ and $y=h(z)^{1/2}$.  
Then we have the same realization of the complex hyperbolic space 
$SU(1,n)/SU(n)$ as in Venkov \cite{venkov}: 
if we write $z_k=x_k+\sqrt{-1}y_k, k=2,...,n$, 
the Riemannian metric is written as 
\begin{equation*}
\dd s^2 = \frac{1}{y^2}\dd y^2 + 
\frac{1}{y^2} \sum_{k=2}^n (\dd x_k^2 + \dd y_k^2) + 
\frac{1}{y^4} \biggl( \dd x_1 + \sum_{k=2}^n 
(x_k \dd y_k - y_k \dd x_k) \biggr)^2, 
\end{equation*}
and the distance function $d(z,z')$ is given by 
\begin{equation*}
(\cosh (d(z,z')))^2 = 
\frac{((y')^2+\Phi)^2 + 4 \varphi^2}{4y^2(y')^2},
\end{equation*}
where 
\begin{equation} \label{ce:phidef} 
\Phi = y^2 + |\wt{z}'-\wt{z}|^2 \qquad \text{\rm and} \qquad 
\varphi=x_1'-x_1+\sum_{k=2}^n(y'_kx_k-x'_ky_k).
\end{equation}
The Laplace-Beltrami operator is written by 
\begin{equation} \label{ce:lap} \begin{split}
\Delta_c = y^4 \frac{\pa^2}{\pa x_1^2} & + y^2 \frac{\pa^2}{\pa y^2} 
- (2n-1) y\frac{\pa}{\pa y} \\ 
 & + y^2 \sum_{k=2}^n \biggl\{ 
\biggl(\frac{\pa}{\pa x_k} + y_k \frac{\pa}{\pa x_1}\biggr)^2 + 
\biggl(\frac{\pa}{\pa y_k} - x_k \frac{\pa}{\pa x_1}\biggr)^2 
\biggr\}. 
\end{split} \end{equation}
\indent Letting $(W^{(2n)},\Bo^{(2n)}, P^{(2n)})$ be 
the $(2n)$-dimensional standard Wiener space and 
denoting an element of $W^{(2n)}$ by 
\begin{equation*}
(B(\cdot),w_2(\cdot),w_3(\cdot),...,w_{2n}(\cdot))\qquad\text{\rm or} 
\qquad (B(\cdot),w_2(\cdot),\widetilde{w}(\cdot)), 
\end{equation*}
we can check that 
the Brownian motion $\{Z(t)\}$ on $\Hy_c^n$ 
with $Z(0)=(x_1,y,z_2,$ $...,$ $z_n)$, 
$z_k$ being identified with $(x_k,y_k)$, is given by 
\begin{equation} \label{ce:sol} \begin{split} 
 & X(t) = x_1 + \int_0^t Y(s)^2 \dd w_2(s) + 2 \sum_{k=2}^n S_k(t) \\
 & Y(t) = y \exp(B(t)-nt) \\
 & X_k(t) = x_k + \int_0^t Y(s) \dd w_{2k-1}(s), \\
 & Y_k(t) = y_k + \int_0^t Y(s) \dd w_{2k}(s), \qquad \qquad k=2,...,n,
\end{split} \end{equation}
where we have used the trivial notations $X(t), Y(t), X_k(t), Y_k(t)$ 
for the components of $Z(t)$ and 
$S_k(t)$ is the stochastic area enclosed 
by $\{(X_k(s),Y_k(s))\}_{0\leqq s \leqq t}$ and its chord, 
\begin{equation*}
S_k(t) = \frac12 \int_0^t (Y_k(s) \dd X_k(s) - X_k(s) \dd Y_k(s)).
\end{equation*}
\indent $\{Y(t)\}$ is again a usual geometric Brownian motion 
with negative drift and $Z(t)$ converges as $t\to\infty$.   
Hence we easily obtain the following central limit theorem.

\medskip

\begin{thm} \label{ct;clt}
For the Brownian motion $\{Z(t)\}$ 
on the $n$-dimensional complex hyperbolic space{\rm ,} 
the probability law of 
$\sqrt{t}^{-1}(d(Z(t),Z(0))-nt)$ converges weakly as $t\to\infty$ 
to the standard normal distribution. 
\end{thm}

\medskip

Next we compute the limiting distribution 
of $(X(t),\widetilde{Z}(t))$ as $t\to\infty$, where 
\begin{equation*}
\wt{Z}(t) = (X_2(t), Y_2(t), ..., X_n(t), Y_n(t)).
\end{equation*}
If we consider the ball model, we obtain the Poisson kernels 
as the densities of the image measures of the uniform measure 
on the sphere by the isometries.   
However, since the same strategy works 
in the more complicated case of the quaternionic hyperbolic space 
whose geometry has not been well understood (see a recent work 
by Kim-Parker \cite{kp} and references cited therein), 
we give the following straightforward computations.   

For this purpose we first fix $t$ and 
consider the characteristic function.  
As in the previous section, we set 
\begin{equation*}
A_t^{(-\mu)}=\int_0^t e^{2B^{(-\mu)}_s} \dd s, \qquad 
\widetilde{A}_t^{(-\mu)}=\int_0^t e^{4B^{(-\mu)}_s} \dd s, 
\end{equation*}
$B^{(-\mu)}_s=B(s)-\mu s$ and $\mu=n$.   
For the stochastic analysis on $\Hy_c^n$ and $\Hy_q^n$, 
we need to consider these two exponential functionals.   
Then, by the expression \eqref{ce:sol}, it is easy to see 
that, for fixed $t>0$, $(X(t),\wt{Z}(t))$ is identical in law with 
\begin{equation*}
\bigl( x_1 + y^2w_2\bigl(\wt{A}_t^{(-\mu)}\bigr) + 
y \phi\bigl(A_t^{(-\mu)}\bigr) 
+ 2y^2\sum\wt{S}_k\bigl(A_t^{(-\mu)}\bigr), 
\widetilde{z}+y\wt{w}\bigl(A_t^{(-\mu)}\bigr) \bigr),
\end{equation*}
where $\sum$ denotes the sum over $k=2,...,n$, 
$\phi(t)=\sum(y_k w_{2k-1}(t)-x_k w_{2k}(t))$ and 
\begin{equation*}
\wt{S}_k(t) = \frac12 \int_0^t (w_{2k}(s) \dd w_{2k-1}(s) - 
w_{2k-1}(s) \dd w_{2k}(s)).
\end{equation*}
Hence we may write, for any bounded continuous function $g$ 
on $\R^{2(n-1)}$, 
\begin{multline*}
E[e^{\sqrt{-1}p X(t)} g(\wt{Z}(t))] 
= E\biggl[e^{\sqrt{-1}p(x+y^2w_2(\wt{A}_t^{(-\mu)})+
y\phi(A_t^{(-\mu)}))} 
g\bigl(\wt{z}+y\wt{w}(A_t^{(-\mu)})\bigr) \\ 
\times E\biggl[\prod_{k=2}^n e^{2\sqrt{-1}py^2\wt{S}_k(A_t^{(-\mu)})} 
\bigg| \{B(s)\}, \wt{w}\bigl(A_t^{(\mu)}\bigr) \biggr] \biggr].
\end{multline*}
Then, applying the L\'evy formula for the characteristic function 
of the stochastic area (cf. \cite{iw}, p.473), we get 
\begin{multline*}
E[e^{\sqrt{-1}p X(t)} g(\wt{z}+\wt{Z}(t))] 
= E\biggl[e^{\sqrt{-1}p(x+y^2w_2(\wt{A}_t^{(-\mu)})+
y\phi(A_t^{(-\mu)}))} 
g\bigl(\wt{z}+y\wt{w}\bigl(A_t^{(-\mu)}\bigr)\bigr) \\
 \times \biggl( \frac{py^2A_t^{(-\mu)}}{\sinh(py^2A_t^{(-\mu)})}
\biggr)^{n-1} 
\exp\biggl(\bigl(1-py^2A_t^{(-\mu)}\coth(py^2A_t^{(-\mu)})\bigr)
\frac{|\wt{w}(A_t^{(-\mu)})|^2}{2A_t^{(-\mu)}}\biggr) \biggr].
\end{multline*}
Moreover we take the conditional expectation 
given $\{B(s)\}_{s\geqq0}$ or $\{Y(s)\}_{s\geqq0}$ to obtain 
\begin{multline*}
E[e^{\sqrt{-1}p X(t)} g(\wt{Z}(t))] 
 = e^{\sqrt{-1}px} E\biggl[ e^{-p^2y^4\wt{A}_t^{(-\mu)}/2} 
\int_{\R^{2(n-1)}} 
e^{\sqrt{-1}p\sum(y_k\xi_k-x_k\eta_k)} \\ 
 \times g(\wt{z}+\zeta) 
 \biggl(\frac{p}{2\pi\sinh(py^2A_t^{(-\mu)})}\biggr)^{n-1} 
 e^{-p\coth(py^2A_t^{(-\mu)})|\zeta|^2/2} 
\dd\mathbf{\xi}\dd\mathbf{\eta} 
\biggr], 
\end{multline*}
where $\zeta=(\mathbf{\xi},\mathbf{\eta})=
(\xi_2,\eta_2,...,\xi_n,\eta_n)$.  

Now we put, for $\mathbf{q}=(q_2,...,q_n), 
\mathbf{r}=(r_2,...,r_n)\in\R^{n-1}$, 
\begin{equation*}
g(\zeta) = \exp(\sqrt{-1}
(\langle \mathbf{q},\mathbf{\xi} \rangle + 
\langle \mathbf{r},\mathbf{\eta} \rangle)).
\end{equation*}
Then, carring out the Gaussian integral with respect to $\mathbf{\xi}$ 
and $\mathbf{\eta}$, we get 
\begin{multline*} 
E\biggl[\exp\biggl\{ 
\sqrt{-1}(pX(t)+ \sum(q_kX_k(t)+r_kY_k(t)))\biggr\}\biggr]\\
 = e^{\sqrt{-1}f} E\biggl[ e^{-p^2y^4\wt{A}_t^{(-\mu)}/2} 
\biggl(\frac{1}{\cosh(py^2A_t^{(-\mu)})}\biggr)^{n-1} 
e^{-F\tanh(py^2A_t^{(-\mu)})} \biggr], 
\end{multline*}
where 
\begin{align*}
 & f=f(p,\mathbf{q},\mathbf{r})= px + \sum(q_kx_k+r_ky_k) \\
\intertext{and}
 & F=F(p,\mathbf{q},\mathbf{r})=
 \sum\frac{(q_k+py_k)^2+(r_k-px_k)^2}{2p}.
\end{align*}
Now, letting $t\to\infty$, we obtain the following.

\medskip

\begin{prop}
For any $p\in\R, \mathbf{q},\mathbf{r}\in\R^{n-1},$ 
one has 
\begin{multline} \label{ce:limit} 
 \lim_{t\to\infty}
 E\biggl[\exp\biggl\{ 
\sqrt{-1}(pX(t) + \sum(q_kX_k(t)+r_kY_k(t)))\biggr\}\biggr] \\ 
 = e^{\sqrt{-1}f} E\biggl[ e^{-p^2y^4\wt{A}^{(-n)}_\infty/2} 
\biggl(\frac{1}{\cosh(py^2A^{(-n)}_\infty)}\biggr)^{n-1} 
e^{-F\tanh(py^2A^{(-n)}_\infty)} \biggr].
\end{multline} 
\end{prop}

\medskip

Denote the right hand side of \eqref{ce:limit} 
by $I(p,\mathbf{q},\mathbf{r})$.  
By using the joint Laplace transform of $A_\infty^{(-n)}$ and 
$\wt{A}_\infty^{(-n)}$ given by Corollary \ref{pc:cor} 
in the appendix, we obtain 
\begin{align*}
 & I(p,\mathbf{q},\mathbf{r}) \\
 & = e^{\sqrt{-1}f} \int_0^\infty 
 \biggl[ e^{-p^2y^4\wt{A}^{(-n)}_\infty/2} 
 \bigg| A^{(-n)}_\infty = u \biggr] 
 \biggl(\frac{1}{\cosh(py^2u)}\biggr)^{n-1} \\
 & \hspace{6cm} \times 
 e^{-F \tanh(py^2u)} P(A^{(-n)}_\infty \in \dd u) \\
 & = e^{\sqrt{-1}f} \int_0^\infty \frac{1}{2^n\Gamma(n)} 
 \biggl(\frac{py^2}{\sinh(py^2u)}\biggr)^{n+1} 
 e^{-py^2\coth(py^2u)/2} \\
 & \hspace{6cm} \times 
 \biggl(\frac{1}{\cosh(py^2u)}\biggr)^{n-1} e^{-F \tanh(py^2u)} \dd u.
\end{align*}
Then, changing the variable, 
we see that $I(p,\mathbf{q},\mathbf{r})$ is equal to 
\begin{multline*} 
 e^{\sqrt{-1}f} \frac{(py^2)^n}{2^n\Gamma(n)} \int_0^\infty 
 \biggl(\frac{1}{\sinh(u)}\biggr)^{n+1}\biggl(\frac{1}{\cosh(u)}\biggr)^{n-1}
 e^{-py^2\coth(u)/2} \\
 \times \exp\biggl(- \sum 
 \frac{(q_k+py_k)^2+(r_k-px_k)^2}{2p} \tanh(u) \biggr)\; \dd u
\end{multline*}
if $p>0$ and to 
\begin{multline*} 
 e^{\sqrt{-1}f} \frac{(-py^2)^n}{2^n\Gamma(n)} \int_0^\infty 
 \biggl(\frac{1}{\sinh(u)}\biggr)^{n+1}\biggl(\frac{1}{\cosh(u)}\biggr)^{n-1}
 e^{py^2\coth(u)/2} \\
 \times \exp\biggl( \sum 
 \frac{(q_k+py_k)^2+(r_k-px_k)^2}{2p} \tanh(u) \biggr)\; \dd u
\end{multline*}
if $p<0$.  
From these expressions, we can take the Fourier inversion 
\begin{equation*}
f_n(x',\wt{z}';z) \equiv \frac{1}{(2\pi)^{2n-1}} \int_{\R^{2n-1}} 
I(p,\mathbf{q},\mathbf{r}) e^{-\sqrt{-1}(px'+\sum(q_kx'_k+r_ky'_k))} 
\dd p \dd q_2 \cdots \dd r_n.
\end{equation*}
\indent For the integral with respect to $q_k$ when $p>0$, 
we note as usual 
\begin{multline*} 
 - \frac{(q_k+py_k)^2}{2p}\tanh(u) + \sqrt{-1}q_k(x_k-x'_k) \\ 
 = - \frac{\tanh(u)}{2p} \bigl( q_k + py_k - 
\sqrt{-1}p(x_k-x'_k)\coth(u)\bigr)^2 \\ 
 \qquad - \sqrt{-1}py_k(x_k-x'_k) - 
\frac{p}{2}(x_k-x'_k)^2 \coth(u).
\end{multline*}
We do the same computations 
also for the other variables and for $p<0$.  
Then, after some manipulations, we obtain 
\begin{align*}
 & f_n(x',\wt{z}';z) \\
 & = \frac{y^{2n}}{(4\pi)^n\Gamma(n)} \int_{\R} |p|^{2n-1} 
 e^{\sqrt{-1}\varphi p} \dd p \int_0^\infty 
 \biggl(\frac{1}{\sinh(u)}\biggr)^{2n} e^{-\Phi|p|\coth(u)/2} \dd u \\
 & = \frac{2y^{2n}}{\pi^n\Gamma(n)} \int_0^\infty p^{2n-1} 
 \cos(2\varphi p) 
 \dd p \int_0^\infty \biggl(\frac{1}{\sinh(u)}\biggr)^{2n} 
 e^{-\Phi p \coth(u)} \dd u,
\end{align*}
where we have made a simple change of variable for the second equality 
and $\varphi$ and $\Phi$ are given by \eqref{ce:phidef}.  

For the second integral, 
we change the variable by $k=\coth(u)$ to obtain 
\begin{equation*}
f_n(x',\wt{z}';z) = \frac{2y^{2n}}{\pi^n\Gamma(n)} \int_0^\infty 
p^{2n-1} \cos(2\varphi p) \dd p \int_1^\infty e^{-\Phi pk} 
(k^2-1)^{n-1} \dd k.
\end{equation*}
Now we recall the following integral representation 
of the modified Bessel function (cf. Lebedev \cite{leb} p.119 or 
\cite{g-r} p.322) 
\begin{equation} \label{ce:bessel} 
K_\nu(z) = 
\frac{\sqrt{\pi}}{\Gamma(\nu+1/2)} \biggl(\frac{z}{2}\biggr)^\nu 
\int_1^\infty e^{-zt} (t^2-1)^{\nu-1/2} \dd t, \qquad \nu>0.
\end{equation}
Then we obtain 
\begin{equation*}
f_n(x',\wt{z}';z) = \frac{2^{n+1/2}y^{2n}}{\pi^{n+1/2}\Phi^{n-1/2}} 
\int_0^\infty p^{n-1/2} \cos(2\varphi p) K_{n-1/2}(\Phi p) \dd p
\end{equation*}
For the integral on the right hand side, 
we may apply the formulae 
\begin{multline*} 
\int_0^\infty x^\lambda K_\mu(ax) \cos(bx) \dd x = 
2^{\lambda-1} a^{-\lambda-1} 
\Gamma\bigl(\frac{\mu+\lambda+1}{2}\bigr) 
\Gamma\bigl(\frac{1+\lambda-\mu}{2}\bigr) \\
 \times F\bigl(\frac{\mu+\lambda+1}{2},\frac{1+\lambda-\mu}{2};
\frac{1}{2};-\frac{b^2}{a^2}\bigr).
\end{multline*}
(cf. \cite{g-r} p.747) and $F(n,a,a;z)=(1-z)^{-n}.$  
Then we obtain 
\begin{equation} \label{ce:poisson} 
f_n(x',\wt{z}';z) = \frac{2^{2n-1}\Gamma(n)y^{2n}}{\pi^n\Phi^{2n}} 
\!\; _1F_0\bigl(n;-\frac{4\varphi^2}{\Phi^2}\bigr) = 
\frac{2^{2n-1}\Gamma(n)y^{2n}}{\pi^n(4\varphi^2+\Phi^2)^n}. 
\end{equation}

\medskip

\begin{thm}[cf. \cite{deb1}] \label{ct:poisson} 
For any $z\in\Hy_c^n,$ $(X(t),\wt{Z}(t))$ converges almost surely 
as $t\to\infty$ 
and the density of the limit distribution on $\R^{2n-1}$ is 
the Poisson kernel given by \eqref{ce:poisson}.
\end{thm}
\section{Quaternionic hyperbolic spaces}

For the quaternion hyperbolic space $Sp(1,n)/(Sp(1)\times Sp(n))$, 
$n\geqq 2$, we follow the conventions in Venkov \cite{venkov}.   
See also Helgason \cite{H}, Lohou\'e-Rychner \cite{LR}, 
and Kim-Parker \cite{kp} for the basic properties.   
For $n\geqq2$, let $\Hy_q^n$ be the upper half space in $\C^{2n}$, 
\begin{equation*}
\Hy_q^n = \{ z=(z_1,z_2,...,z_{2n})=(z_1,\widetilde{z})\in\C^{2n}; 
\mathrm{Im}(z_1)>0\}, 
\end{equation*}
with the Riemannian metric 
\begin{multline*}
ds^2 = \frac{\dd y^2}{y^2} + 
\frac{1}{y^2} \sum_{k=2}^n(\dd z_k \dd \overline{z_k} + 
\dd z_{n+k} \dd \overline{z_{n+k}}) \\
+ \frac{1}{y^4} \biggl( \dd x_1 + \mathrm{Im} \sum_{k=2}^n 
(\overline{z_k} \dd z_k + \overline{z_{n+k}} \dd z_{n+k})\biggl)^2 \\
+ \frac{1}{y^4} \bigg| \dd z_{n+1} + \sum_{k=2}^n 
(z_{n+k} \dd z_k - z_k \dd z_{n+k}) \bigg|^2,
\end{multline*}
where $z_1=x_1+\sqrt{-1}y$.  
We will write $z_k=x_k+\sqrt{-1}y_k$ for $k=2,...,2n$.  
Note that the first and $(n+1)$-th components, 
$z_1$ and $z_{n+1}$ play special roles.   

The volume element is $y^{-4n-3}\dd x_1 \dd y \prod_{k=2}^{2n} 
\dd x_k \dd y_k$ and 
the distance function $d(z,z')$ is given by 
\begin{equation} \label{qe:dist} 
\bigl( \cosh(d(z,z')) \bigr)^2 = 
\frac{((y')^2+\Phi)^2+4(\varphi_1^2+\varphi_2^2+\varphi_3^2)}
{4y^2(y')^2},
\end{equation}
where 
\begin{equation} \label{qe:phi} \begin{split}
 & \Phi= y^2 + \sum_{k=2}^n(|z'_k-z_k|^2 + 
|z'_{n+k}-z_{n+k}|^2), \\
 & \varphi_1 = x_1'-x_1 + \sum_{k=2}^n \bigl( (y'_k x_k - x'_k y_k) + 
 (y'_{n+k}x_{n+k}-x'_{n+k}y_{n+k}) \bigr), \\
& \varphi_2 = x'_{n+1}-x_{n+1} + 
\sum_{k=2}^n \bigl( (x'_k x_{n+k} - x'_{n+k} x_k) + 
 (y'_{n+k} y_{k}- y'_{k} y_{n+k}) \bigr), \\
& \varphi_3 = y'_{n+1}-y_{n+1} + 
\sum_{k=2}^n \bigl( (x'_k y_{n+k} - y'_{n+k} x_k) + 
 (y'_k x_{n+k} - x'_{n+k} y_k) \bigr). 
\end{split} \end{equation}
Note that $\varphi_i$'s do not depend on $y$.  

The Laplace-Beltrami operator $\Delta_q$ may be written 
in a convenient way as 
\begin{equation} \label{qe:lap} \begin{split}
 \Delta_q = & y^4 \frac{\pa^2}{\pa x_1^2} + y^2 \frac{\pa^2}{\pa y^2} 
 - (4n+1) y \frac{\pa}{\pa y} + 
y^4 \biggl( \frac{\pa^2}{\pa x_{n+1}^2} + \frac{\pa^2}{\pa y_{n+1}^2} 
\biggr) \\ 
 & \qquad + y^2 \sum_{k=2}^n \biggl[ \biggl( \frac{\pa}{\pa x_k} + 
y^k \frac{\pa}{\pa x_1} - x_{n+k} \frac{\pa}{\pa x_{n+1}} 
- y_{n+k} \frac{\pa}{\pa y_{n+1}} \biggr)^2 \\
 & \qquad \qquad \qquad + \biggl( \frac{\pa}{\pa y_k} - 
x^k \frac{\pa}{\pa x_1} + y_{n+k} \frac{\pa}{\pa x_{n+1}} 
- x_{n+k} \frac{\pa}{\pa y_{n+1}} \biggr)^2 \\
 & \qquad \qquad \qquad + \biggl( \frac{\pa}{\pa x_{n+k}} + 
y_{n+k} \frac{\pa}{\pa x_1} + x_{k} \frac{\pa}{\pa x_{n+1}} 
+ y_{k} \frac{\pa}{\pa y_{n+1}} \biggr)^2 \\ 
 & \qquad \qquad \qquad + \biggl( \frac{\pa}{\pa y_{n+k}} - 
x_{n+k} \frac{\pa}{\pa x_1} - y_{k} \frac{\pa}{\pa x_{n+1}} 
+ x_{k} \frac{\pa}{\pa y_{n+1}} \biggr)^2 \biggr]. 
\end{split} \end{equation}
Note that the coefficients of $\pa^2/\pa x_1\pa x_{n+1}, 
\pa^2/\pa x_1 \pa y_{n+1},$ $\pa^2/\pa x_{n+1} \pa y_{n+1}$ are 
zero.  
We can describe the same story as for the complex hyperbolic space 
if we consider a $4\times 4$ skew-symmetric matrix 
instead of $2$-simensional one.   

At first we give an explcit expression for the Brownian motion, 
the diffusion process with generator $\Delta_q/2$, on $\Hy_q^n$.  
Let $(W^{(4n)},\Bo^{(4n)}, P^{(4n)})$ be the $(4n)$-dimensional 
Wiener space and denote an element in $W^{(4n)}$ by 
\begin{equation*} \begin{split} 
(B_1(\cdot),B(\cdot) & ,w_{2,1}(\cdot),w_{2,2}(\cdot),...,
w_{n,1}(\cdot),w_{n,2}(\cdot), \\ 
 & B_2(\cdot),B_3(\cdot),w_{n+2,1}(\cdot),w_{n+2,2}(\cdot),..., 
w_{2n,1}(\cdot),w_{2n,2}(\cdot)). 
\end{split} \end{equation*}
Then we can check that 
the Brownian motion $(X(t),Y(t),\wt{Z}(t))$ 
starting from $(x_1,y,\wt{z})$ is given by 
\begin{align*}
 & X_1(t) = x_1 + \int_0^t Y(s)^2 \dd B_1(s) \\
 & \qquad \qquad + \sum_{k=2}^n \int_0^t \{ 
 Y_k(s) \dd X_k(s) - X_{k}(s) \dd Y_k(s) \\
 & \qquad \qquad \qquad \qquad \qquad 
 + Y_{n+k}(s) \dd X_{n+k}(s) 
 - X_{n+k}(s) \dd Y_{n+k}(s) \} , \\
 & Y(t) = y \exp(B(t)-(2n+1)t), \\
 & X_k(t) = x_k + \int_0^t Y(s) \dd w_{k,1}(s) , \\
 & Y_k(t) = y_k + \int_0^t Y(s) \dd w_{k,2}(s) , 
 \qquad \quad k=2,...,n, \\
 & X_{n+1}(t) = x_{n+1} + \int_0^t Y(s)^2 \dd B_2(s) \\
 & \qquad \qquad + \sum_{k=2}^n \int_0^t \{ 
 - X_{n+k}(s) \dd X_k(s) + X_{k}(s) \dd X_{n+k}(s) \\
 & \qquad \qquad \qquad \qquad \qquad 
 + Y_{n+k}(s) \dd Y_{k}(s) 
 - Y_{k}(s) \dd Y_{n+k}(s) \} , \\
 & Y_{n+1}(t) = y_{n+1} + \int_0^t Y(s)^2 \dd B_3(s) \\
 & \qquad \qquad + \sum_{k=2}^n \int_0^t \{ 
 - Y_{n+k}(s) \dd X_k(s) + X_{k}(s) \dd Y_{n+k}(s) \\
 & \qquad \qquad \qquad \qquad \qquad 
 - X_{n+k}(s) \dd Y_{k}(s) 
 + Y_{k}(s) \dd X_{n+k}(s) \}\\
 & X_{n+k}(t) = x_{n+k} + \int_0^t Y(s) \dd w_{n+k,1}(s) , \\
 & Y_{n+k}(t) = y_{n+k} + \int_0^t Y(s) \dd w_{n+k,2}(s) , 
 \qquad \quad k=2,...,n. 
\end{align*}
\indent Then, from \eqref{qe:dist}, it is easy to show 
the following central limit theorem.  

\medskip

\begin{thm} \label{qt:clt} 
The probability law of $(d(Z(t),Z(0))-(2n+1)t)/\sqrt{t}$ converges 
weakly as $t\to\infty$ to the standard normal distribution.  
\end{thm}

\medskip

Next we show that $(X(t),\wt{Z}(t))$ converges in law as $t\to\infty$. 
To identify the limit distribution, we set 
\begin{equation*}
f_n(x',\wt{z}';z) = \frac{2^{4n+1}\Gamma(2n)}{\pi^{2n}} 
\frac{y^{2(2n+1)}}
{(\Phi^2+4(\varphi_1^2+\varphi_2^2+\varphi_3^2))^{2n+1}} ,
\end{equation*}
where $\Phi$ and $\varphi_i$'s are given by \eqref{qe:phi}.  
$f_n$ is the Poisson kernel of 
the boundary $\pa\Hy_q^n=\{y=0\}$.  

\medskip

\begin{thm} \label{qt:poisson} 
$(X_1(t),\wt{Z}(t)),$ valued in $\R\times\C^{2(n-1)},$ 
converges almost surely as $t\to\infty$ and 
the density of the limit distribution is 
given by $f_n(x',\wt{z}';z).$ 
\end{thm}

\medskip

In the following we give a proof of Theorem \ref{qt:poisson}.   
At first we consider the characteristic function 
of $(X_1(t),\wt{Z}(t))$ for fixed $t$.  
For convenience we put 
\begin{equation*}
X_k^0(t)=X_k(t)-x_k=\int_0^t Y(s)\dd w_{k,1}(s), \ 
Y_k^0(t)=Y_k(t)-y_k=\int_0^t Y(s)\dd w_{k,2}(s),
\end{equation*}
and 
\begin{equation*}
\theta_k = 
\begin{pmatrix} x_k \\ y_k \\ x_{n+k} \\ y_{n+k} \end{pmatrix}, \quad 
\Theta_k(t) = 
\begin{pmatrix} X_k(t) \\ Y_k(t) \\ 
X_{n+k}(t) \\ Y_{n+k}(t) \end{pmatrix}, \quad 
\Theta_k^0(t) = 
\begin{pmatrix} X^0_k(t) \\ Y^0_k \\ 
X^0_{n+k} \\ Y^0_{n+k}(t) \end{pmatrix}. 
\end{equation*}
Moreover, $\xi=\!\;^t(\xi_1,\xi_2,\xi_3)\in\R^3, 
w_k=\!\;^t(u_k,v_k,u_{n+k},v_{n+k})\in\R^4, 
w=(w_2,...,w_n),$ we set 
\begin{align*}
 & \Psi(t) = \xi_1 X_1(t) + \xi_2 X_{n+1}(t) + \xi_3 Y_{n+1}(t), \\
 & U_k(t) = \langle w_k, \Theta_k(t) \rangle, \qquad 
 U^0_k(t) = \langle w_k, \Theta^0_k(t) \rangle.
\end{align*}
We throughout denote by $^tQ$ the transpose of a matrix $Q$.  
Then the characteristic function is 
\begin{align*}
\varphi(t) & = E\bigl[ 
\exp\bigl\{ \sqrt{-1} 
\bigl( \xi_1 X_1(t) + \xi_2 X_{n+1}(t) + \xi_3 Y_{n+1}(t) \bigr) \\
 & \qquad \qquad + \sqrt{-1} \sum \bigl(u_k X_k(t) + v_k Y_k(t) 
 + u_{n+k} X_{n+k}(t) + v_{n+k} Y_{n+k}(t) \bigr) \bigr\} \bigr] \\
 & = E\bigl[ \exp\bigl( \sqrt{-1} ( \Psi(t) + \sum U_k(t) \bigr) 
\bigr],
\end{align*}
where the summation is taken over $k=2,...,n$.  

To compute the characteristic function, 
we introduce a $4\times 4$ skew symmetric matrix $\Xi$ given by 
\begin{equation*}
\Xi = \begin{pmatrix} 
0      & \xi_1  & -\xi_2 & -\xi_3 \\
-\xi_1 & 0      & -\xi_3 & \xi_2 \\
\xi_2  & \xi_3  & 0      & \xi_1 \\
\xi_3  & -\xi_2 & -\xi_1 & 0         \end{pmatrix}.
\end{equation*}
Then we have 
\begin{equation*} \begin{split}
\Psi(t) + \sum U_k(t) = \psi & + \int_0^t Y(s)^2 
(\xi_1 \dd B_1(s) + \xi_2 \dd B_2(s) + \xi_3 \dd B_3(s)) \\ 
 & + \sum \langle \Xi \theta_k + w_k, \Theta^0_k(t) \rangle + 
 \sum \int_0^t \langle \Xi \Theta^0_k(s), 
 \dd \Theta^0_k(s) \rangle,
\end{split} \end{equation*}
where $\psi=\xi_1 x_1 + \xi_2 x_{n+1} + \xi_3 y_{n+1} + 
\sum\langle w_k,\theta_k \rangle$.  
Note that $\{\sum_{j=1}^3\xi_j B_j(s)\}$ is identical in law with 
$\{|\xi|B_1(s)\}$, $|\xi|=(\xi_1^2+\xi_2^2+\xi_3^2)^{1/2}$.   

The eigenvalues of $\Xi$ are $\pm\sqrt{-1}|\xi|$ and 
the multiplicities are two.   
Moreover there exists an orthogonal matrix $Q$ 
such that $^tQ\Xi Q=K$ is of the standard form.   
We take 
\begin{equation*}
Q = 
\begin{pmatrix} 
0 & \frac{\xi_1}{|\xi|} & \frac{\sqrt{\xi_2^2+\xi_3^2}}{|\xi|} & 0 \\
1 & 0 & 0 & 0 \\
0 & \frac{\xi_3}{|\xi|} 
  & \frac{-\xi_1\xi_3}{\sqrt{\xi_2^2+\xi_3^2}|\xi|} 
  & \frac{\xi_2}{\sqrt{\xi_2^2+\xi_3^2}} &  \\
0 & \frac{-\xi_2}{|\xi|} 
  & \frac{\xi_1\xi_2}{\sqrt{\xi_2^2+\xi_3^2}|\xi|} 
  & \frac{\xi_3}{\sqrt{\xi_2^2+\xi_3^2}}
\end{pmatrix} 
\quad \text{\rm and} \quad 
K = 
\begin{pmatrix} 0     & -|\xi| & 0 & 0 \\
                |\xi| &  0     & 0 & 0 \\
                0     &  0     & 0 & -|\xi| \\
                0     &  0     & |\xi| & 0 \end{pmatrix}.
\end{equation*}
We also put $\wh{w}_k=\!\;^tQw_k, \wh{\theta}_k=\!\;^tQ\theta_k$ and 
\begin{equation*}
\wh{\Theta}^0_k(t) = (\wh{X}_k^0(t),\wh{Y}^0_k(t),\wh{X}^0_{n+k}(t),
\wh{Y}^0_{n+k}(t)) = \!\;^tQ \wh{\Theta}^0_k(t).
\end{equation*}
By the rotation invariance of the probability law of Brownian motions, 
we see that $\{\wh{\Theta}_k^0(s)\}$ is a simple time change 
of a $4$-dimensional standard Brownian motion.  

Under these notations, we have 
\begin{equation*}
\langle \Xi \theta_k + w_k, \Theta^0_k(t) \rangle = 
\langle K\wh{\theta}_k + \wh{w}_k, \wh{\Theta}^0_k(t) \rangle
\end{equation*}
and 
\begin{equation} \label{qe:area} \begin{split} 
 \int_0^t \langle \Xi \Theta^0_k(s),\dd \Theta^0_k(s) \rangle & = 
\int_0^t \langle K \wh{\Theta}_k^0(s), \dd\wh{\Theta}_k^0(s) \rangle \\
 & = |\xi|\int_0^t\{ \wh{X}_k^0(s)\dd\wh{Y}_k^0(s) - 
 \wh{Y}^0_k(s)\dd\wh{X}^0_k(s) \\
 & \qquad \qquad + \wh{X}^0_{n+k}(s)\dd\wh{Y}^0_{n+k}(s) 
 - \wh{Y}^0_{n+k}(s)\dd\wh{X}^0_{n+k}(s) \}.
\end{split} \end{equation}
As in the previous sections, we set 
\begin{equation*}
A_t^{(-\mu)} = \int_0^t \exp(2B_s^{(-\mu)}) \dd s 
\qquad \text{\rm and} \qquad 
\wt{A}_t^{(-\mu)} = \int_0^t \exp(4B_s^{(-\mu)}(s)) \dd s, 
\end{equation*}
$B_s^{(-\mu)}=B(s)-\mu s$ and $\mu=2n+1$.  
Then, by taking the conditional expectation 
given $\{Y(s)\}$ and $\wh{\Theta}^0_k(t),k=2,...,n$, and 
applying the L\'evy formula, we obtain 
\begin{align*}
 & \varphi(t) = e^{\sqrt{-1}\psi} E\biggl[ \exp \biggl( 
-\frac12 |\xi|^2 y^4 \wt{A}_t^{(-\mu)} + 
\sqrt{-1}\sum\langle K\wh{\theta}_k+\wh{w}_k,\wh{\Theta}_k(t) 
\rangle \biggr) \\
 & \qquad \qquad \times \prod_{k=2}^n 
 E\biggl[ \exp\biggl( \sqrt{-1} \int_0^t 
\langle \Xi \Theta^0_k(s),\dd \Theta^0_k(s) \rangle \biggr) \bigg| 
\{Y(s)\}, \wh{\Theta}_k(t) \biggr] \biggr] \\
 & = e^{\sqrt{-1}\psi} E\biggl[ \exp \biggl( 
-\frac12 |\xi|^2 y^4\wt{A}_t^{(-\mu)} + 
\sqrt{-1}\sum\langle K\wh{\theta}_k+\wh{w}_k,y\wh{W}_k(A_t^{(-\mu)}) 
\rangle \biggr) \\
 & \times 
 \biggl(\frac{|\xi|y^2A_t^{(-\mu)}}{\sinh(|\xi|y^2A_t^{(-\mu)})}
\biggr)^{2(n-1)} \\
 & \times \exp\biggl( 
\bigl( 1-|\xi|y^2A_t^{(-\mu)}\coth(|\xi|y^2A_t^{(-\mu)}) \bigr)
\frac{|\wh{W}(A_t^{(-\mu)})|^2}{2A_t^{(-\mu)}} \biggr) \biggr] \\
 & = e^{\sqrt{-1}\psi} 
E\biggl[ \exp \biggl( -\frac12 |\xi|^2 y^4\wt{A}_t^{(-\mu)}\biggr) 
\int_{\R^{4(n-1)}} 
\biggl(\frac{|\xi|}{2\pi\sinh(|\xi|y^2A_t^{(-\mu)})}
\biggr)^{2(n-1)} \\
 & \times 
\exp\biggl( \sqrt{-1} \sum \biggl( 
\langle K\wh{\theta}_k+\wh{w}_k,\zeta_k \rangle 
- \frac{1}{2} |\xi| \coth(|\xi|y^2A_t^{(-\mu)})|\zeta_k|^2 \biggr) 
\biggr) \ \dd\zeta_2\cdots\dd\zeta_n \biggr],
\end{align*}
where $\{(\wh{W}_2(s),...,\wh{W}_{n}(s))\}$ is 
a $4(n-1)$-dimensional standard Brownian motion, 
independent of $\{Y(s)\}$ or $\{B(s)\}$.  
We carry out the Gaussian integral over $\R^{4(n-1)}$ to obtain 
\begin{equation*} \begin{split}
 & \varphi(t) \\ 
 & = e^{\sqrt{-1}\psi} E\biggl[ 
\biggl( \frac{1}{\cosh(|\xi|y^2A_t^{(-\mu)})} \biggr)^{2(n-1)} 
e^{-|\xi|^2y^4\wt{A}_t^{(-\mu)}/2 - 
F\tanh(|\xi|y^2A_t^{(-\mu)})/2|\xi|} \biggr] ,
\end{split} \end{equation*}
where $F=\sum|K\wh{\theta}_k+\wh{w}_k|^2$.  
Hence, letting $t$ tend to $\infty$, we obtain 
\begin{equation*} \begin{split}
 & \lim_{t\to\infty}\varphi(t) \\ 
 & = e^{\sqrt{-1}\psi} E\biggl[ 
\biggl( \frac{1}{\cosh(|\xi|y^2A^{(-\mu)}_\infty)} \biggr)^{2(n-1)} 
e^{-|\xi|^2y^4\wt{A}^{(-\mu)}_\infty/2 - 
F\tanh(|\xi|y^2A^{(-\mu)}_\infty)/2|\xi|} \biggr].
\end{split} \end{equation*}
\indent Now, applying \eqref{pe:cor}, 
we obtain the explicit expression 
for the Fourier transform of the Poisson kernel $f_n$.  

\medskip

\begin{prop} \label{qp:fou}
Under the notations above{\rm,} 
the Fourier transform of the limit distribution 
of $(X(t),\wt{Z}(t))$ as $t\to\infty$, is given by 
\begin{multline*} 
\Phi(\xi,w) = 
\frac{e^{\sqrt{-1}\psi}(|\xi|y^2)^{2n+1}}{2^{2n+1}\Gamma(2n+1)} 
\int_0^\infty \biggl(\frac{1}{\cosh(u)}\biggr)^{2(n-1)} 
\biggl(\frac{1}{\sinh(u)}\biggr)^{2(n+1)} \\
 \times 
 \exp\biggl(-\frac12 |\xi|y^2\coth(u)-\frac{F}{2|\xi|}\tanh(u)\biggr) 
 \dd u.
\end{multline*}
\end{prop} 

\medskip

We invert the Fourier transform.   
That is, setting $x'=(x'_1,x'_{n+1},y'_{n+1})$, 
$\theta'_k=\!\;^t(x'_k, y'_k, x'_{n+k}, y'_{n+k})$ and 
\begin{equation*}
\psi' = \xi_1 x'_1 + \xi_2 x'_{n+1} + \xi_3 y'_{n+1} + 
\sum \langle w_k,\theta'_k \rangle, 
\end{equation*}
we compute 
\begin{equation*} \begin{split} 
 & \bar{f}_n(x',\wt{z}';z) \\
 & \equiv \frac{1}{(2\pi)^{4n-1}} \int_{\R^{4n-1}} 
 e^{\sqrt{-1}(\psi-\psi')} \dd\xi \dd w \\ 
 & \qquad \qquad \quad \qquad \times \int_0^\infty 
 \frac{(|\xi|y^2)^{2n+1}}{2^{2n+1}\Gamma(2n+1)} 
 \biggl(\frac{1}{\cosh(u)}\biggr)^{2(n-1)} 
 \biggl(\frac{1}{\sinh(u)}\biggr)^{2(n+1)} \\
 & \qquad \qquad \qquad \quad \times 
 \exp\biggl( -\frac{1}{2}|\xi|y^2\coth(u) - 
 \frac{\tanh(u)}{2|\xi|}\sum|K\wh{\theta}_k\wh{w}_k|^2\biggr)\ \dd u.
\end{split} \end{equation*}
Recall the definitions, 
$\wh{w}_k=\!\;^tQ w_k$ and $\wh{\theta}_k=\!\;^tQ \theta_k$.  
Then, changing the order of the integrations, we have 
\begin{equation*} \begin{split} 
 & \bar{f}_n(x',\wt{z}';z) \\
 & = \frac{y^{2(2n+1)}}{(2\pi)^{4n-1}2^{2n+1}\Gamma(2n+1)} 
 \int_0^\infty \biggl(\frac{1}{\cosh(u)}\biggr)^{2(n-1)} 
 \biggl(\frac{1}{\sinh(u)}\biggr)^{2(n+1)} \dd u \\
 & \qquad \times \int_{\R^3} e^{\sqrt{-1}\langle \xi,x-x' \rangle} 
e^{-|\xi|y^2\coth(u)/2} |\xi|^{2n+1} \dd\xi \\
 & \qquad \times \int_{\R^{4(n-1)}} 
e^{\sqrt{-1}\sum\langle \wh{w}_k,\wh{\theta}_k-\wh{\theta}'_k\rangle
-\tanh(u)\sum|K\wh{\theta}_k+\wh{w}_k|^2/2|\xi|}\prod_{k=2}^n \dd w_k.
\end{split} \end{equation*}
We can easily carry out the third Gaussian integral 
since $Q\in O(4)$ and we obtain 
\begin{multline*}
\bar{f}_n(x',\wt{z}';z) 
= \frac{y^{2(2n+1)}}{(4\pi)^{2n+1}\Gamma(2n+1)} 
\int_0^\infty \biggl(\frac{1}{\sinh(u)}\biggr)^{4n} \dd u \\
 \times \int_{\R^3} e^{\sqrt{-1}\{ \langle \xi,x-x' \rangle + 
|\xi|\phi(\wh{\theta},\wh{\theta}') \} 
-|\xi| \Phi \coth(u)/2} |\xi|^{4n-1} \dd\xi,
\end{multline*}
where $\Phi=y^2+\sum|\wh{\theta}_k-\wh{\theta}'_k|^2=
y^2+\sum|\theta_k-\theta'_k|^2$ and 
\begin{equation*}
\phi(\wh{\theta},\wh{\theta}') = 
\langle K \wh{\theta}_k, \wh{\theta}_k' \rangle = 
\sum(\wh{y}'_k\wh{x}_k-\wh{x}'_k\wh{y}_k 
+ \wh{y}'_{n+k}\wh{x}_{n+k}-\wh{x}'_{n+k}\wh{y}_{n+k}).
\end{equation*}
\indent For the right hand side, 
changing the variables by $k=\coth(u)$, 
we first compute the integral in $u$.  
Then, by using formula \eqref{ce:bessel} again, we get 
\begin{equation*} \begin{split} 
 \int_0^\infty \biggl(\frac{1}{\sinh(u)}\biggr)^{4n} 
 e^{-|\xi| \Phi \coth(u)/2} \dd u 
 & = \int_1^\infty (k^2-1)^{2n-1} e^{-|\xi| \Phi k/2} \dd k \\
 & = \frac{\Gamma(2n)}{\sqrt{\pi}} 
\biggl(\frac{4}{|\xi|\Phi}\biggr)^{2n-1/2} 
K_{2n-1/2}\biggl(\frac{|\xi|\phi}{2}\biggr).
\end{split} \end{equation*}
Moreover, by definitions, we see 
\begin{equation*}
\langle \xi,x'-x \rangle + |\xi| \phi(\wh{\theta},\wh{\theta}') = 
-\langle \xi,\varphi\rangle,
\end{equation*}
where $\varphi=\!\:^t(\varphi_1,\varphi_2,\varphi_3)$ is 
given by \eqref{qe:phi}. 

Combining these identities, we obtain 
\begin{align*}
 & \bar{f}_n(x',\wt{z}';z) \\
 & = \frac{y^{2(2n+1)}}{8(2n+1)\pi^{2n+3/2}\Phi^{2n-1/2}} \int_{\R^3} 
e^{-\sqrt{-1}\langle \xi,\varphi \rangle} 
K_{2n-1/2}\biggl(\frac{\Phi}{2}|\xi|\biggr) 
|\xi|^{2n-1/2} \dd\xi \\ 
 & = \frac{y^{2(2n+1)}}{8(2n+1)\pi^{2n+3/2}\Phi^{2n-1/2}} \int_{\R^3} 
e^{\sqrt{-1}|\varphi|\xi_3} 
K_{2n-1/2}\biggl(\frac{\Phi}{2}|\xi|\biggr) 
|\xi|^{2n-1/2} \dd\xi.
\end{align*}
Moreover, changing the variables by the spherical coordinate, 
we obtain 
\begin{equation*} \begin{split} 
\bar{f}_n(x',\wt{z}';z) & = 
\frac{4\pi y^{2(2n+1)}}{8(2n+1)\pi^{2n+3/2}\Phi^{2n-1/2}|\varphi|} \\
 & \qquad \qquad \times 
\int_0^\infty r^{2n+1/2} K_{2n-1/2}\biggl(\frac{\Phi}{2}r\biggr) 
\sin(|\varphi|r) \dd r.
\end{split} \end{equation*}
For the integral on the right hand side, the following formula 
is available (cf. \cite{g-r}, p.747):
\begin{multline*} 
\int_0^\infty x^\lambda K_\mu(ax) \sin(bx) \dd x = 
2^{\lambda} a^{-\lambda-2} b \times
\Gamma\bigl(\frac{2+\mu+\lambda}{2}\bigr) 
\Gamma\bigl(\frac{2+\lambda-\mu}{2}\bigr) \\
 \times F\bigl(\frac{2+\mu+\lambda}{2},\frac{2+\lambda-\mu}{2};
\frac{3}{2};-\frac{b^2}{a^2}\bigr).
\end{multline*}
In our case $\lambda=2n+1/2, \mu=2n-1/2$ and 
\begin{equation*}
F\bigl(\frac{2+\mu+\lambda}{2},\frac{2+\lambda-\mu}{2};
\frac{3}{2};-\frac{b^2}{a^2}\bigr) = 
\biggl(1+\frac{b^2}{a^2}\biggr)^{-(2n+1)}. 
\end{equation*}
Hence we may apply this identity and we arrive at our result 
\begin{equation*}
\bar{f}_n(x,\wt{z}';z) = \frac{2^{4n+1}\Gamma(2n)}{\pi^{2n}} 
\frac{y^{2(2n+1)}}{(\Phi^2+4|\varphi|^2)^{2n+1}}.
\end{equation*}
\appendix
\section{Perpetual integrals of geometric Brownian motion}

In this appendix we consider two perpetual integrals of 
geometric Brownian motions.   
Let $B=\{B(t)\}_{t\geqq0}$ be a one-dimensional Brownian motion 
with $B_0=0$ defined on a probability space $(\Omega, \F,P)$.  
For $\mu>0$, we set $B^{(-\mu)}=\{B_t^{(-\mu)}\equiv B(t)-\mu t\}$, 
a Brownian motion with negative constant drift $-\mu$.  
Then Dufresne's perpetual integral is defined by 
\begin{equation} \label{pe:expfnl1} 
A_\infty^{(-\mu)} = \int_0^\infty \exp(2B_s^{(-\mu)}) \dd s. 
\end{equation}
We also consider another integral 
\begin{equation*}
a_\infty^{(-\mu)} = \int_0^\infty \exp(B_s^{(-\mu)}) \dd s. 
\end{equation*}

Then the following is known:

\medskip

\begin{thm}[Dufresne \cite{duf}] \label{pt:marg}
Let $\gamma_\mu$ be a gamma random variable 
whose density is given by $(\Gamma(\mu))^{-1}x^{\mu-1}e^{-x}.$  
Then $A_\infty^{(-\mu)}$ is distribued as $(2\gamma_\mu)^{-1}$ 
and{\rm ,} accordingly{\rm ,} $a_\infty^{(-\mu)} \eil 
2(\gamma_{2\mu})^{-1}.$
\end{thm}

\medskip

\begin{rem}
Several different proofs of this theorem are known.   
In particular, see Yor \cite{yor-sur}.  
The density of the exponential functional 
$A_t=\int_0^t\exp(2B_s) ds$ for fixed $t$ has been obtained 
by Yor \cite{yor-some} and 
the joint distribution of $(A_t,a_t)$ in an obvious notation 
has been studied in \cite{ams}.   
See also \cite{survey1, survey2, yor-coll} 
for several results and applications 
of these perpetual integrals and exponential functionals.   
Recently Baudoin-O'Connell \cite{bo} has shown several formulae,
including \eqref{pe:joint} below, 
for the exponential type Wiener and discussed their close relation 
to the theory of quantum Toda lattice.   
\end{rem}

\medskip

What we need in Sections 4 and 5 is the following explicit 
expression for the conditional Laplace transform 
of $A_\infty^{(-\mu)}$ given $a_\infty^{(-\mu)}$, 
which was originally obtained by Yor \cite{yor-two}.   
We set 
\begin{equation*}
f_1(v) = \frac{2^{2\mu}}{\Gamma(2\mu)} v^{-(2\mu+1)} e^{-2/v}, 
\qquad v>0, 
\end{equation*}
which is the density of the random variable $a_\infty^{(-\mu)}$ or 
$2/\gamma_{2\mu}$.   

\medskip

\begin{thm} \label{pt:joint}
For $\lambda>0$ and $v>0,$ it holds that 
\begin{equation} \label{pe:joint} \begin{split}
E\biggl[ & \exp\biggl( - \frac{1}{2}\lambda^2 A_\infty^{(-\mu)} 
\biggr) \bigg| a_{\infty}^{(-\mu)}=v\biggr] f_1(v) \\ 
 & = \frac{1}{2\Gamma(2\mu)} 
\biggl(\frac{\lambda}{\sinh(\lambda v/2)}\biggr)^{2\mu+1} 
\exp\biggl(-\lambda \coth\biggl(\frac{\lambda v}{2}\biggr) \biggr).
\end{split} \end{equation}
\end{thm}

\medskip

We have this nice result only for the particular choice of 
$A_t$ and $a_t$, that is, 
it is available only when the ratio of the coefficients 
in the exponential functionals is two.  

We give another proof of the theorem for completeness. 
Note that by letting $\lambda$ tend to $0$ in \eqref{pe:joint}, 
we obtain Theorem \ref{pt:marg}.  

For this purpose, 
we consider the Brownian motion $\{B_t^{(\mu)}=B_t+\mu t\}$ 
with the opposite positive drift and set $X_x(s)=x\exp(B_s^{(\mu)})$, 
which defines a diffusion process with infinitesimal generator 
\begin{equation*}
\frac12 x^2 \frac{\dd^2}{\dd x^2} + \biggl(\frac{1}{2}+\mu\biggr) 
x\frac{\dd}{\dd x}.
\end{equation*}
Letting $\tau_z$ be the first hitting time of $\{X_x(s)\}$ 
at $z$, we set for $\lambda>0$ and $\kappa\in\R$ 
\begin{equation*}
v_z(x) = E\biggl[ \exp\biggl( -\frac{\lambda^2}{2} 
\int_0^{\tau_z}X_x(s)^{-2} \dd s + \lambda\kappa \int_0^{\tau_z} 
X_x(s)^{-1}\dd s\biggl)\biggl].
\end{equation*}
\indent In \cite{survey1} we have considered the case of $\kappa=0$ 
and showed that $v_z(x)$ may be represented by means of 
the modified Bessel function 
to give another proof of Theorem \ref{pt:marg}.   
Following the same line, we first show a representation 
for $v_z(x)$ by means of the Whittaker function.  

Let $W_{\kappa,\mu}$ be a Whittaker function: 
if $\mu-\kappa+1/2>0$, 
\begin{equation} \label{pe:int-rep}
W_{\kappa,\mu}(z) = 
\frac{e^{-z/2}z^{\mu+1/2}}{\Gamma(\mu-\kappa+1/2)} 
\int_0^\infty e^{-zt} t^{\mu-\kappa-1/2} (1+t)^{\mu+\kappa-1/2} \dd t.
\end{equation}
From this expression it is easy to see 
$\lim_{z\to\infty}W_{\kappa,\mu}(z)=0$ when $|\kappa|$ is small.  
We also recall that $W_{\kappa,\mu}$ solves the equation 
\begin{equation*}
W''(z) + \biggl(-\frac{1}{4} + \frac{\kappa}{z} - 
\frac{\mu^2-(1/4)}{z^2}\biggr) W(z) = 0.
\end{equation*}

\medskip

\begin{prop} \label{pp:whittaker}
For $\mu>0,\lambda>0$ and $\kappa\in\R,$ it holds that 
\begin{equation} \label{pe:whittaker}
v_z(x) = \biggl(\frac{z}{x}\biggr)^{\mu-1/2} 
\frac{W_{\kappa,\mu}(2\lambda/x)}{W_{\kappa,\mu}(2\lambda/z)}.
\end{equation}
\end{prop}

\noindent{\it Proof}.  
We have only to consider the case of $\kappa<0$.   
The general case can be shown from this case 
by analytic continuation in $\kappa$.   
Note that, if $\kappa<0$, $v_z(x)$ is monotone increasing 
in $x (>z)$.   

At first we note that $v_z(x)$ is a solution for 
\begin{equation*}
\frac12 x^2 v''(x) + \biggl(\frac{1}{2}+\mu\biggr) x v'(x) = 
\biggl(\frac{\lambda^2}{2x^2}-\frac{\lambda\kappa}{x}\biggr) v(x)
\end{equation*}
and satisfies 
\begin{equation} \label{pe:bdry}
v_z(x)\bigg|_{x=z} = 1 \qquad \text{\rm and} \qquad 
\lim_{x\downarrow0}v_z(x)=0.
\end{equation}
We now change the variable by $\xi=\lambda/x$ and set 
\begin{equation*}
v_z(x)=\xi^{\mu-1/2} \phi(\xi).
\end{equation*}
Then, by straightforward computations, 
we see that $\phi$ satisfies 
\begin{equation*}
\phi''(\xi) + \biggl( -1 + \frac{2\kappa}{\xi} - 
\frac{\mu^2-(1/4)}{\xi^2}\biggr) \phi(\xi) = 0.
\end{equation*}
By considering the boundary conditions \eqref{pe:bdry}, 
we can easily show 
\begin{equation*}
\phi(\xi) = \biggl(\frac{z}{\lambda}\biggr)^{\mu-1} 
\frac{W_{\kappa,\mu}(2\xi)}{W_{\kappa,\mu}(2\lambda/z)}
\end{equation*}
and hence the result \eqref{pe:whittaker}. \hfill $\square$

\medskip

\begin{prop} \label{pp:lt}
For $\mu>0,\lambda>0$ and $\kappa\in\R,$ it holds that 
\begin{equation} \label{pe:lt}
E\biggl[ \exp\biggl( -\frac12 \lambda^2 A_\infty^{(-\mu)} 
+\lambda\kappa a_\infty^{(-\mu)}\biggr) \biggr] = 
\frac{\Gamma(\mu-\kappa+1/2)}{\Gamma(2\mu)} 
(2\lambda)^{\mu-1/2} W_{\kappa,\mu}(2\lambda).
\end{equation}
\end{prop}

\noindent{\it Proof}.  
By the symmetry of the probability law of Brownian motion, 
$\{-B_t\}\eil\{B_t\}$, we have 
\begin{align*}
\lim_{z\to\infty}v_z(1) & = 
E\biggl[ \exp\biggl( -\frac12 \lambda^2 \int_0^\infty 
e^{-2B_s^{(\mu)}} \dd s + \lambda \kappa 
\int_0^\infty e^{-B_s^{(\mu)}} \dd s \biggr) \biggr] \\
 & = E\biggl[ \exp\biggl( -\frac12 \lambda^2 \int_0^\infty 
e^{2B_s^{(-\mu)}} \dd s + \lambda \kappa 
\int_0^\infty e^{B_s^{(-\mu)}} \dd s \biggr) \biggr].
\end{align*}
On the other hand, by using the fact on the Whittaker function 
\begin{equation} \label{pe:asym} 
W_{\kappa,\mu}(z) = \frac{\Gamma(2\mu)}{\Gamma(\mu-\kappa+1/2)} 
z^{-\mu+1/2} (1+o(1)) \quad \text{\rm as} \quad 
z\downarrow 0,
\end{equation}
we see from the expression \eqref{pe:whittaker}
\begin{equation*}
\lim_{z\to\infty}v_z(1) = \frac{\Gamma(\mu-\kappa+1/2)}{\Gamma(2\mu)} 
(2\lambda)^{\mu-1/2} W_{\kappa,\mu}(2\lambda). 
\end{equation*}
\hfill $\square$ 

\begin{rem}
The asymptotic behavior \eqref{pe:asym} of $W_{\kappa,\mu}$ can 
be easily shown by the definition of the Whittaker functions.
\end{rem}

Now we are in a position to complete our proof of 
Theorem \ref{pt:joint}.  
By \eqref{pe:int-rep} and \eqref{pe:lt}, we have 
\begin{equation*} \begin{split} 
E\biggl[ \exp \biggl(-\frac12 \lambda^2 A_\infty^{(-\mu)} & + 
\lambda \kappa a_\infty^{(-\mu)} \biggr) \biggr] \\ 
 & = \frac{e^{-\lambda}}{\Gamma(2\mu)} (2\lambda)^{2\mu} 
\int_0^\infty e^{-2\lambda t} t^{\mu-\kappa-1/2} 
(1+t)^{\mu+\kappa-1/2} \dd t.
\end{split} \end{equation*}
Now we change the variable by $e^{\lambda v}=1+t^{-1}$ or 
$t=(e^{\lambda v}-1)^{-1}$.  
Then some elementary computations show that this integral is equal to 
\begin{equation*}
\frac{\lambda^{2\mu+1}}{2\Gamma(2\mu)} \int_0^\infty 
e^{-\lambda\coth(\lambda v/2)} 
\biggl(\frac{1}{\sinh(\lambda v/2)}\biggr)^{2\mu+1} 
e^{\lambda\kappa v} \dd v.
\end{equation*}
This completes our proof since 
\begin{equation*} \begin{split}
E\biggl[ \exp \biggl(-\frac12 \lambda^2 A_\infty^{(-\mu)} & + 
\lambda \kappa a_\infty^{(-\mu)} \biggr) \biggr] \\ 
 & = \int_0^\infty 
E\biggl[\exp\biggl(-\frac12 \lambda^2 A_\infty^{(-\mu)}\biggr)\bigg|
a_\infty^{(-\mu)}=v \biggr] f_1(v) e^{\lambda\kappa v} \dd v.
\end{split} \end{equation*}
\hfill $\square$

\medskip

\begin{cor} \label{pc:cor}
Define another perpetual integral $\widetilde{A}_\infty^{(-\mu)}$ by 
\begin{equation} \label{pe:expfnl3} 
\widetilde{A}_\infty^{(-\mu)} = \int_0^\infty 
\exp(4B_s^{(-\mu)}) \dd s 
\end{equation}
and let $f_2(v)$ be the density of $A_\infty^{(-\mu)}$ or 
$(2\gamma_\mu)^{-1}.$  Then one has 
\begin{equation} \label{pe:cor} \begin{split}
E\biggl[ \exp
\biggl(-\frac{1}{2}\lambda^2 \widetilde{A}_\infty^{(-\mu)} \biggr) 
\bigg| & A_{\infty}^{(-\mu)}=v\biggr] f_2(v) \\ 
 & = \frac{1}{2^\mu\Gamma(\mu)} 
\biggl(\frac{\lambda}{\sinh(\lambda v)}\biggr)^{\mu+1} 
\exp\biggl(-\frac{\lambda}{2} \coth(\lambda v) \biggr).
\end{split} \end{equation}
\end{cor}

\bigskip

\begin{flushleft}
Graduate School of Information Science \\
Nagoya University \\
Chikusa-ku, Nagoya 464-8601, Japan \\
E-mail:matsu@is.nagoya-u.ac.jp 
\end{flushleft}

\end{document}